\documentclass[a4paper,12pt]{article}
\usepackage{amssymb,amsmath}

\newcommand{\N}{\mathbb{N}}

\newcommand{\twosum}[2]{\sum_{\substack{#1\\#2}}}

\newcommand{\ep}{\varepsilon}
\newtheorem{theorem}{Theorem}
\newtheorem{lemma}{Lemma}
\renewcommand{\mod}[1]{\hspace{-2.9mm}\pmod{#1}}

\begin{document}
\title{Fractional Moments of Dirichlet $L$-Functions}
\author{D.R. Heath-Brown\\Mathematical Institute, Oxford}
\date{}
\maketitle
\section{Introduction}

Mean-values of the type
\[I_k(T):=\int_0^T|\zeta(\tfrac12+it)|^{2k}dt,\]
with positive non-integral values of $k$, have been investigated by a 
number of authors, including Ramachandra \cite{Ram},
\cite{RamFenn}, Conrey and Ghosh \cite{CG} and Heath-Brown
\cite{frac1}.  In particular the above papers by Ramachandra show,
under the Riemann Hypothesis, that
\[I_k(T)\gg_k T(\log T)^{k^2}\;\;\; (T\ge 2)\]
for all real $k\ge 0$, and that
\[I_k(T)\ll_k T(\log T)^{k^2}\;\;\; (T\ge 2)\]
for all real $k\in[0,2]$.

It is natural to ask about the corresponding problem for Dirichlet
$L$-functions in $q$-aspect, that is to say to investigate
\[M_k(q):=\twosum{\chi\mod{q}}{\chi\not=\chi_0}|L(\tfrac12,\chi)|^{2k}\]
for positive real $k$.  However rather little is known about this in general.
The method of Rudnick and Soundararajan \cite{RS1}, enables
one to show unconditionally that
\[M_k(q)\gg_k \phi(q)(\log q)^{k^2}\]
for rational $k\ge 1$, at least when $q$ is prime. 
It is annoying that the range $0\le k <1$ is not covered by
this approach.

The present paper will prove results in the reverse direction,
motivated by the author's work \cite{frac1}.  We establish the
following theorems.
\begin{theorem}\label{T1}
Assuming the Generalized Riemann Hypothesis we have
\[M_k(q)\ll_k \phi(q)(\log q)^{k^2}\]
for all $k\in(0,2)$.
\end{theorem}
\begin{theorem}\label{T2}
Unconditionally we have
\[M_k(q)\ll_k \phi(q)(\log q)^{k^2}\]
for any $k$ of the form $k=1/v$, with $v\in\N$.
\end{theorem}
Thus taking $v=2$ we have
\[\sum_{\chi\mod{q}}|L(\tfrac12,\chi)|\ll \phi(q)(\log q)^{1/4}\]
in particular.

The approach in \cite{frac1} is based on a convexity theorem for
mean-value integrals, which appears to have no analogue for character
sums.  We therefore work with integrals, and extract the sum $M_k(q)$
at the end.  While we can give lower bounds for the integrals that
occur, as well as upper bounds, it is not clear how to give a lower 
bound for $M_k(q)$ in terms of an integral.

This work arose from a number of conversations with Dr H.M. Bui, and
would not have been undertaken without his prompting.
It is a pleasure to acknowledge his contribution.

\section{Mean-Value Integrals}

Throughout our argument we will write $v=1$ for the proof of Theorem 
\ref{T1}, and $v=k^{-1}$ in handling Theorem \ref{T2}.  In both cases
the primary mean-value integral we will work with is
\[J(\sigma,\chi):=\int_{-\infty}^{\infty}|L(\sigma+it,\chi)|^{2k}
|W(\sigma+it)|^6dt,\]
where the weight function $W(s)$ is defined by
\[W(s):=\frac{q^{\delta(s-1/2)}-1}{(s-1/2)\log q},\]
with $\delta>0$ to be specified later, see (\ref{dc1}) and (\ref{dc2}).  
We emphasize that, for the rest of this
paper, all constants implied by the Vinogradov $\ll$ symbol will be
uniform in $\sigma$ for the ranges specified.  However they will be 
allowed to depend on the values of $k$ and $\delta$, so that the
symbol $\ll$ should be read as $\ll_{k,\delta}$ throughout.  

In addition to the integral $J(\sigma,\chi)$ we will use
\[K(\sigma,\chi):=\int_{-\infty}^{\infty}|S(\sigma+it,\chi)|^2
|W(\sigma+it)|^6dt,\]
where 
\[S(s):=\sum_{n\le q}d_k(n)\chi(n)n^{-s} \] 
Notice here that a little
care is needed in defining $d_k(n)$ when $k$ is not an integer, see
\cite[\S 2]{frac1}.

When $\chi$ is a non-principal character the function $L(s,\chi)$ is
entire.  Moreover, if we assume the Generalized Riemann Hypothesis
then there are no zeros for $\sigma>\tfrac12$, so that one can define
a holomorphic extension of 
\[L(s,\chi)^k=\sum_{m=1}^{\infty}d_k(m)\chi(m)m^{-s}\;\;\;(\sigma>1)\] 
in the half-plane $\sigma>\tfrac12$.  Having defined $L(s,\chi)^k$ in
this way we now set
\[G(\sigma,\chi):=
\int_{-\infty}^{\infty}|L(\sigma+it,\chi)^{k}-S(\sigma+it,\chi)|^2
|W(\sigma+it)|^6dt, \;\;\;(\sigma>\tfrac12).\]
This integral will be used in the proof of Theorem \ref{T1}, while for
the unconditional Theorem \ref{T2} we will employ
\[H(\sigma,\chi):=
\int_{-\infty}^{\infty}|L(\sigma+it,\chi)-S(\sigma+it,\chi)^v|^{2/v}
|W(\sigma+it)|^6dt.\]

In addition to $J(\sigma,\chi), K(\sigma,\chi), G(\sigma,\chi)$ and
$H(\sigma,\chi)$ we will consider their averages over non-principal
characters,
\[J(\sigma):=\twosum{\chi\mod{q}}{\chi\not=\chi_0}J(\sigma,\chi),\;\;\;\;\;
K(\sigma):=\twosum{\chi\mod{q}}{\chi\not=\chi_0}K(\sigma,\chi)\]
\[G(\sigma):=\twosum{\chi\mod{q}}{\chi\not=\chi_0}G(\sigma,\chi),\;\;\;\;\;
\mbox{and}\;\;\;\;\;
H(\sigma):=\twosum{\chi\mod{q}}{\chi\not=\chi_0}H(\sigma,\chi).\]

To derive estimates relating values of these integrals we begin with
the following convexity estimate of Gabriel \cite[Theorem 2]{gab}.
\begin{lemma}\label{con}
Let $F$ be a complex-valued function which is regular in the strip
$\alpha<\Re(z)<\beta$, and continuous for 
$\alpha\le\Re(z)\le\beta$.  Suppose that $|F(z)|$ tends to
zero as $|\Im(z)|\rightarrow\infty$, uniformly for 
$\alpha\le\Re(z)\le\beta$.
Then for any $\gamma\in[\alpha,\beta]$ and any $a>0$ we have
\[I(\gamma)\le I(\alpha)^{(\beta-\gamma)/(\beta-\alpha)}
I(\beta)^{(\gamma-\alpha)/(\beta-\alpha)}\]
where
\[I(\eta):=\int_{-\infty}^{\infty}|F(\eta+it)|^adt.\]
\end{lemma}
The inequality should be interpreted appropriately if any of the integrals
diverge. From Lemma \ref{con} we will deduce the following variant.
\begin{lemma}\label{newcon}
Let $f$ and $g$ be complex-valued functions which are regular in the strip
$\alpha<\Re(z)<\beta$, and continuous for 
$\alpha\le\Re(z)\le\beta$.  Let $b$ and $c$ be positive real numbers.
Suppose that $|f(z)|^b|g(z)|^c$ and $|g(z)|$ tend
to zero as
$|\Im(z)|\rightarrow\infty$, uniformly for $\alpha\le\Re(z)\le\beta$.
Set
\[I(\eta):=\int_{-\infty}^{\infty}|f(\eta+it)|^b|g(\eta+it)|^cdt.\]
Then for any $\gamma\in[\alpha,\beta]$ we have
\begin{equation}\label{in1}
I(\gamma)\le I(\alpha)^{(\beta-\gamma)/(\beta-\alpha)}
I(\beta)^{(\gamma-\alpha)/(\beta-\alpha)}.
\end{equation}
\end{lemma}
To deduce Lemma \ref{newcon} from Lemma \ref{con} we choose a rational
number $p/q>c/b$, and apply Lemma \ref{con} with $F=f^qg^p$ and
$a=b/q$.  Since 
\[|F|=(|f|^b|g|^c)^{q/b}|g|^{p-cq/b}\]
with $p-cq/b>0$,
we deduce that $|F|$ tends to zero as
$|\Im(z)|\rightarrow\infty$, uniformly for $\alpha\le\Re(z)\le\beta$.
We then obtain an inequality of the same shape as (\ref{in1}), but
with the exponent $c$ replaced by $bp/q$.  Lemma \ref{newcon} then
follows on choosing a sequence of rationals $p_n/q_n$ tending
downwards to $c/b$.

We now apply Lemma \ref{newcon} to $J(\sigma,\chi)$.  When $\sigma=3/2$ we have
\[W(s)\ll q^{\delta}/(1+|t|)\]
whence we trivially obtain
\[J(\tfrac32,\chi)\ll q^{6\delta}.\]
An immediate application of Lemma \ref{newcon} therefore yields
\[J(\sigma,\chi)\ll
J(\tfrac12,\chi)^{3/2-\sigma}q^{6\delta(\sigma-1/2)}\]
for $\tfrac12\le\sigma\le\tfrac32$, whence we trivially deduce that
\[J(\sigma)\ll J(\tfrac12)^{3/2-\sigma}q^{6\delta(\sigma-1/2)},\]
by H\"{o}lder's inequality.  Since 
\begin{equation}\label{Gfin}
J^f\le \left(\frac{\log q}{q}\right)^{1-f}\left(\frac{q}{\log q}+J\right)
\ll q^{-(1-\delta)(1-f)}\left(\frac{q}{\log q}+J\right)
\end{equation}
for any $J\ge 0$ and any $f\in[0,1]$, we conclude as follows.
\begin{lemma}\label{Jconvlem}
We have
\[J(\sigma)\ll q^{-(1-7\delta)(\sigma-1/2)}
\left(\frac{q}{\log q}+J(\tfrac12)\right)\]
for $\tfrac12\le\sigma\le \tfrac32$.
\end{lemma}

To obtain a second estimate involving $J(\sigma,\chi)$ we use Lemma
\ref{newcon} to show that if
$\tfrac12\le\sigma\le\tfrac34$ and $1-\sigma\le\gamma\le\sigma$ then
\[J(\gamma,\chi)\le J(\sigma,\chi)^{(\gamma-1+\sigma)/(2\sigma-1)}
J(1-\sigma,\chi)^{(\sigma-\gamma)/(2\sigma-1)}.\]
An application of H\"{o}lder's inequality then shows that
\[J(\gamma)\le J(\sigma)^{(\gamma-1+\sigma)/(2\sigma-1)}
J(1-\sigma)^{(\sigma-\gamma)/(2\sigma-1)}.\]
To handle $J(1-\sigma,\chi)$ we will use
the functional equation for $L(s,\chi)$.  If $\psi$ is primitive, with
conductor $q_1$, this yields
\[L(1-\sigma+it,\psi)\ll
(1+|t|)^{\sigma-1/2}q_1^{\sigma-1/2}|L(\sigma+it,\psi)|\]
for $\tfrac12\le\sigma\le\tfrac34$ say.  Thus if $\psi$ induces a
character $\chi$ modulo $q$ we will have
\[L(1-\sigma+it,\chi)\ll
(1+|t|)^{\sigma-1/2}q_1^{\sigma-1/2}\rho|L(\sigma+it,\chi)|\]
with
\[\rho=\prod_{p\mid q_2}
\left(\frac{|1-\chi(p)p^{-\sigma-it}|}{|1-\chi(p)p^{\sigma-1-it}|}\right),\]
where $q_2=q/q_1$.  Thus
\[\log \rho\le 
(2\sigma-1)\sum_{p\mid q_2}\frac{\log p}{p^{1-\sigma}-1}.\]
However
\[\sum_{p|m}\frac{\log p}{p^{1/4}-1}\le\tfrac12\log m\]
for all sufficiently large $m$, whence $\rho\ll q_2^{\sigma-1/2}$.
We therefore conclude that
\[L(1-\sigma+it,\chi)\ll
(1+|t|)^{\sigma-1/2}q^{\sigma-1/2}|L(\sigma+it,\chi)|\]
when $\tfrac12\le\sigma\le\tfrac34$, for any character 
$\chi$ modulo $q$, whether primitive or not.

We now deduce that
\[J(1-\sigma,\chi)\hspace{11cm}\]
\[\ll q^{2k(\sigma-1/2)}
\int_{-\infty}^{\infty}|L(\sigma+it,\chi)|^{2k}(1+|t|)^{2k(\sigma-1/2)}
|W(1-\sigma+it)|^6dt.\]
The presence of the factor $(1+|t|)^{2k(\sigma-1/2)}$ is
inconvenient.  However, since $0<k<2$ we have
\[(1+|t|)^{2k(\sigma-1/2)}|W(1-\sigma+it)|^6\ll
(\log q)^{-6}|t|^{-2},\]
for $|t|\ge 1$ and $\tfrac12\le\sigma\le 1$.  It follows that
\[J(1-\sigma,\chi)\ll q^{2k(\sigma-1/2)}\left(J(\sigma,\chi)+
(\log q)^{-6}J^*(\sigma,\chi)\right),\]
where
\[J^*(\sigma,\chi):=\int_{-\infty}^{\infty}|L(\sigma+it,\chi)|^{2k}
\frac{dt}{1+t^2}.\]
Thus
\[J(1-\sigma)\ll q^{2k(\sigma-1/2)}\left(J(\sigma)+
(\log q)^{-6}J^*(\sigma)\right)\]
with
\[J^*(\sigma):=\twosum{\chi\mod{q}}{\chi\not=\chi_0}
\int_{-\infty}^{\infty}|L(\sigma+it,\chi)|^{2k}
\frac{dt}{1+t^2}.\]
Finally we observe that
\[J(\sigma)^{(\gamma-1+\sigma)/(2\sigma-1)}\left\{J(\sigma)+
(\log
q)^{-6}J^*(\sigma)\right\}^{(\sigma-\gamma)/(2\sigma-1)}\hspace{3cm}\]
\[\hspace{4cm}\le J(\sigma)+(\log q)^{-6}J^*(\sigma).\]
On comparing our results we therefore conclude that
\begin{equation}\label{in2}
J(\gamma)\ll q^{k(\sigma-\gamma)}\left(J(\sigma)+
 (\log q)^{-6}J^*(\sigma)\right).
\end{equation}

We have now to consider $J^*(\sigma)$.  It was shown by Montgomery
\cite[Theorem 10.1]{mont} that
\[\sum_{\chi\mod{q}}\!\!\!\!\!\raisebox{1ex}{*}\;\;
\int_{-T}^T|L(\tfrac12+it,\chi)|^4dt\ll
\phi(q)T(\log qT)^4\]
for $T\ge 2$, where $\Sigma^*$ indicates that only primitive
characters are to be considered.  (It should be noted that there is a
misprint in the statement of \cite[Theorem 10.1]{mont}, in that
$L(\tfrac12+it,\chi)$ should be replaced by $L(\sigma+it,\chi)$.
However we are only interested in the case $\sigma=\tfrac12$.
moreover, in the proof of \cite[Theorem 10.1]{mont}, at the top of
page 83, the reference to Theorem 6.3 should be to Theorem 6.5.)

If $\chi$ is an imprimitive character modulo $q$, induced by a
primitive character $\psi$ with conductor $q_1$, then
\[|L(\tfrac12+it,\chi)|^4\le |L(\tfrac12+it,\psi)|^4\prod_{p\mid
q,\,p\nmid q_1}(1+p^{-1/2})^4.\]
Thus if $\Sigma^{(1)}$ indicates summation over all characters $\chi$
modulo $q$ for which the conductor has a given value $q_1$, we will have
\[\sum_{\chi}\raisebox{1.5ex}{(1)}\int_{-T}^T|L(\tfrac12+it,\chi)|^4dt\ll
\phi(q_1)T(\log q_1T)^4\prod_{p\mid
q,\,p\nmid q_1}(1+p^{-1/2})^4.\]
If we now sum for $q_1|q$ we obtain
\[\sum_{\chi\mod{q}}\int_{-T}^T|L(\tfrac12+it,\chi)|^4dt\ll
T(\log qT)^4f(q),\]
where
\[f(q)=\sum_{q_1|q}\phi(q)\prod_{p\mid
q,\,p\not\mid q_1}(1+p^{-1/2})^4.\]
The function $f$ is multiplicative, with
\[f(p^e)=(1+p^{-1/2})^4+\phi(p)+\phi(p^2)+\ldots+\phi(p^e)=
p^e\left(1+O(p^{-3/2})\right).\]
Thus $f(q)\ll q$ and we conclude that
\[\sum_{\chi\mod{q}}\int_{-T}^T|L(\tfrac12+it,\chi)|^4dt\ll
qT(\log qT)^4.\]

We may now deduce that if $f(s)=L(s,\chi)^2s^{-1}$ then
\[\sum_{\chi\mod{q}}\int_{-\infty}^{\infty}|f(\tfrac12+it)|^2dt\ll
q(\log q)^4.\]
Moreover the trivial bound $L(s,\chi)\ll 1$ for $\sigma=3/2$ shows that
\[\sum_{\chi\mod{q}}\int_{-\infty}^{\infty}|f(\tfrac32+it)|^2dt\ll q.\]
We can therefore apply Lemma \ref{con}, together with H\"{o}lder's
inequality, to deduce that
\[\sum_{\chi\mod{q}}\int_{-\infty}^{\infty}|f(\sigma+it)|^2dt\ll
q(\log q)^4\]
uniformly for $\tfrac12\le\sigma\le\tfrac32$.  A final application of
H\"{o}lder's inequality then implies that
\[J^*(\sigma)\ll q(\log q)^4.\]

We can now insert this into (\ref{in2}) and deduce as follows.
\begin{lemma}\label{Jconvlem2}
We have
\[J(\gamma)\ll q^{k(\sigma-\gamma)}\left(\frac{q}{\log q}+J(\sigma)\right)\]
for $\tfrac12\le\sigma\le 1$ and $1-\sigma\le\gamma\le\sigma$.
\end{lemma}

We now turn our attention to $G(\sigma,\chi)$ and $H(\sigma,\chi)$.
By Lemma \ref{newcon} we have
\[G(\sigma,\chi)\le G(\tfrac12,\chi)^{3/2-\sigma}G(\tfrac32,\chi)^{\sigma-1/2}
\;\;\;(\tfrac12\le\sigma\le\tfrac32)\]
for non-principal characters $\chi$ modulo $q$.  We then find via 
H\"{o}lder's inequality that
\begin{equation}\label{G1}
G(\sigma)\le
G(\tfrac12)^{3/2-\sigma}G(\tfrac32)^{\sigma-1/2}
\end{equation}
Since
\[W(\tfrac32+it)\ll q^{\delta}(1+|t|)^{-1}\]
we see that
\[G(\tfrac32,\chi)\ll q^{6\delta}
\int_{-\infty}^{\infty}|L(\tfrac32+it,\chi)^{k}-S(\tfrac32+it,\chi)|^2
\frac{dt}{1+|t|^2}.\]
However
\[L(\tfrac32+it,\chi)^{k}-S(\tfrac32+it,\chi)=
\sum_{n> q}d_k(n)\chi(n)n^{-3/2-it} \] 
whence
\[\int_{-\infty}^{\infty}|L(\tfrac32+it,\chi)^{k}-S(\tfrac32+it,\chi)|^2
\frac{dt}{1+|t|^2}\hspace{5cm}\]
\[=\pi\sum_{m,n>q}d_k(m)d_k(n)\chi(m)\overline{\chi(n)}
\min\left(m^{-1/2}n^{-5/2}\,,\,n^{-1/2}m^{-5/2}\right).\]
It follows that
\[\sum_{\chi\mod{q}}
\int_{-\infty}^{\infty}|L(\tfrac32+it,\chi)^{k}-S(\tfrac32+it,\chi)|^2
\frac{dt}{1+|t|^2}\hspace{2cm}\]
\[=\pi\phi(q)\twosum{m,n>q}{q|m-n,\, (mn,q)=1}d_k(m)d_k(n)
\min\left(m^{-1/2}n^{-5/2}\,,\,n^{-1/2}m^{-5/2}\right)
\]
To estimate this double sum we use that fact that
$d_k(n)\ll_{\ep}n^{\ep}$ for any fixed $\ep>0$.  This leads to the bound
\[\twosum{m,n>q}{q|m-n}d_k(m)d_k(n)
\min\left(m^{-1/2}n^{-5/2}\,,\,n^{-1/2}m^{-5/2}\right)
\ll_{\ep}q^{2\ep-2}.\]
It therefore follows that
\[\sum_{\chi\mod{q}}
\int_{-\infty}^{\infty}|L(\tfrac32+it,\chi)^{k}-S(\tfrac32+it,\chi)|^2
\frac{dt}{1+|t|^2}\ll_{\ep}q^{2\ep-1}.\]
Inserting this bound into (\ref{G1}) we obtain
\[G(\sigma)\ll_{\ep}
G(\tfrac12)^{3/2-\sigma}q^{(\sigma-1/2)(6\delta+2\ep-1)}.\]
Using (\ref{Gfin}) again, we see that
\[G(\sigma)\ll_{\ep}
q^{1-2\sigma+(7\delta+2\ep)(\sigma-1/2)}
\left(\frac{q}{\log q}+G(\tfrac12)\right)\]
for $\sigma\in[\tfrac12,\tfrac32]$.  The positive number $\ep$ is at
our disposal, and we choose it to be $\ep=\delta/2$, whence
\[G(\sigma)\ll q^{-(1-4\delta)(2\sigma-1)}
\left(\frac{q}{\log q}+G(\tfrac12)\right).\]

The treatment of $H(\sigma,\chi)$ is similar.  This time, 
since $k=1/v$, we have
\[H(\tfrac32,\chi)\le
\left\{\int_{-\infty}^{\infty}|W(\tfrac32+it)|^6dt\right\}^{1-k}\hspace{3cm}\]
\[\hspace{2cm}\times
\left\{\int_{-\infty}^{\infty}|L(\tfrac32+it,\chi)-S(\tfrac32+it,\chi)^v|^2
|W(\tfrac32+it)|^6dt\right\}^k\]
by H\"{o}lder's inequality.  
The first integral on the right is trivially $O(q^{6\delta})$.  Moreover
\[L(\tfrac32+it,\chi)-S(\tfrac32+it,\chi)^v=
\sum_{n> q}a_k(n)\chi(n)n^{-3/2-it} \] 
with certain coefficients $a_k(n)\ll_{\ep}n^{\ep}$.  The argument
then proceeds as before, noting that
\[\twosum{m,n>q}{q|m-n}a_k(m)a_k(n)
\min\left(m^{-1/2}n^{-5/2}\,,\,n^{-1/2}m^{-5/2}\right)
\ll_{\ep}q^{2\ep-2}.\]
It follows that
\[\sum_{\chi\mod{q}}\int_{-\infty}^{\infty}
|L(\tfrac32+it,\chi)-S(\tfrac32+it,\chi)^v|^2|W(\tfrac32+it)|^6dt\ll 
q^{2\delta+\ep-1}.\]
we then deduce, by the same line of argument as before, that
\[H(\sigma)\ll q^{-(k-4\delta)(2\sigma-1)}
\left(\frac{q}{\log q}+H(\tfrac12)\right)\]
for $\sigma\in[\tfrac12,\tfrac32]$.

We record these results formally in the following lemma.
\begin{lemma}\label{GHconvlem}
For $\sigma\in[\tfrac12,\tfrac32]$ we have
\[G(\sigma)\ll q^{-(1-4\delta)(2\sigma-1)}
\left(\frac{q}{\log q}+G(\tfrac12)\right)\]
and
\[H(\sigma)\ll q^{-(k-4\delta)(2\sigma-1)}
\left(\frac{q}{\log q}+H(\tfrac12)\right).\]
\end{lemma}

We end this section by considering $K(\sigma)$.  We have
\[K(\sigma)\le\sum_{\chi\mod{q}}K(\sigma,\chi)
=\sum_{m,n\le  q}\frac{d_k(m)d_k(n)}{(mn)^{\sigma}}S(m,n)I(m,n),\]
where
\[S(m,n)=\sum_{\chi\mod{q}}\chi(m)\overline{\chi(n)}\]
and
\[I(m,n)=
\int_{-\infty}^{\infty}\left(\frac{n}{m}\right)^{it}|W(\sigma+it)|^6dt.\]
Evaluating the sum $S(m,n)$ we find that
\begin{eqnarray*}
\lefteqn{\sum_{m,n\le
    q}\frac{d_k(m)d_k(n)}{(mn)^{\sigma}}S(m,n)I(m,n)}
\hspace{4cm}\\
&=&\phi(q)\twosum{m,n\le  q}{q\mid m-n,\, (mn,q)=1}
\frac{d_k(m)d_k(n)}{(mn)^{\sigma}}I(m,n)\\
&=&\phi(q)\twosum{n\le  q}{(n,q)=1}\frac{d_k(n)^2}{n^{2\sigma}}
\int_{-\infty}^{\infty}|W(\sigma+it)|^6dt.
\end{eqnarray*}
We then observe that
\[\twosum{n\le  q}{(n,q)=1}\frac{d_k(n)^2}{n^{2\sigma}}\le
\sum_{n\le  q}\frac{d_k(n)^2}{n}\ll (\log q)^{k^2},\]
and that
\[\int_{-\infty}^{\infty}|W(\sigma+it)|^6dt\ll 
q^{3\delta(2\sigma-1)}(\log q)^{-1}.\]
These bounds allow us to conclude as follows.
\begin{lemma}\label{Kbound}
For $\tfrac12\le\sigma\le\tfrac32$ we have
\[K(\sigma)\ll\phi(q)q^{3\delta(2\sigma-1)}(\log q)^{k^2-1}.\]
\end{lemma}

\section{Proof of the Theorems}

By definition of $G(\sigma,\chi)$ and $H(\sigma,\chi)$ we have
\[J(\sigma)\ll K(\sigma)+G(\sigma)\]
under the Generalized Riemann Hypothesis, and
\[J(\sigma)\ll K(\sigma)+H(\sigma)\]
unconditionally.  In view of Lemma \ref{GHconvlem} these produce
\[J(\sigma)\ll K(\sigma)+q^{-(1-4\delta)(2\sigma-1)}
\left(\frac{q}{\log q}+G(\tfrac12)\right)\]
and
\[J(\sigma)\ll K(\sigma)+q^{-(k-4\delta)(2\sigma-1)}
\left(\frac{q}{\log q}+H(\tfrac12)\right)\]
respectively.  However we also have
\[G(\tfrac12)\ll K(\tfrac12)+J(\tfrac12)\]
and
\[H(\tfrac12)\ll K(\tfrac12)+J(\tfrac12)\]
from the definitions again, so that
\[J(\sigma)\ll K(\sigma)+q^{-(1-4\delta)(2\sigma-1)}
\left(\frac{q}{\log q}+K(\tfrac12)+J(\tfrac12)\right)\]
and
\[J(\sigma)\ll K(\sigma)+q^{-(k-4\delta)(2\sigma-1)}
\left(\frac{q}{\log q}+K(\tfrac12)+J(\tfrac12)\right)\]
in the two cases respectively.

If we now call on Lemma \ref{Kbound} we then find that
\begin{eqnarray*}
J(\sigma)&\ll& \phi(q)q^{3\delta(2\sigma-1)}(\log q)^{k^2-1}
+q^{-(1-4\delta)(2\sigma-1)}
\left(\frac{q}{\log q}+J(\tfrac12)\right)\\
&\ll& q^{4\delta(2\sigma-1)}\left(\phi(q)(\log q)^{k^2-1}+
q^{1-2\sigma}J(\tfrac12)\right)
\end{eqnarray*}
under the Generalized Riemann Hypothesis, since
\begin{equation}\label{phiin}
\frac{q}{\log q}\ll \phi(q)(\log q)^{k^2-1}
\end{equation}
for $0<k<2$.  Similarly we have
\[J(\sigma)\ll q^{4\delta(2\sigma-1)}\left(\phi(q)(\log q)^{k^2-1}+
q^{k(1-2\sigma)}J(\tfrac12)\right)\]
unconditionally.

Finally we apply Lemma \ref{Jconvlem2} with
$\gamma=\tfrac12$ and use (\ref{phiin}) again, to deduce that
\[J(\sigma)\ll q^{4\delta(2\sigma-1)}\left(\phi(q)(\log q)^{k^2-1}+
q^{-(2-k)(\sigma-1/2)}J(\sigma)\right)\]
under the Generalized Riemann Hypothesis.  Similarly we may derive the
unconditional bound
\[J(\sigma)\ll q^{4\delta(2\sigma-1)}\left(\phi(q)(\log q)^{k^2-1}+
q^{-k(\sigma-1/2)}J(\sigma)\right).\]
We are now ready to choose our value of $\delta$.  For Theorem \ref{T1}
we take
\begin{equation}\label{dc1}
\delta=\frac{2-k}{10},
\end{equation}
and for Theorem \ref{T2} we choose
\begin{equation}\label{dc2}
\delta=\frac{k}{10}.
\end{equation}
Then in either case we will have
\[J(\sigma)\ll q^{4\delta(2\sigma-1)}\phi(q)(\log q)^{k^2-1}+
q^{-\delta(2\sigma-1)}J(\sigma).\]
We write $c_k$ for the implied constant in this last estimate, and note
that $c_k$ depends only on $k$.  We then take 
\[\sigma=\sigma_0:=\frac{1}{2}+\frac{\kappa}{\log q}\]
with
\[\kappa=(2\delta)^{-1}\max(1\,,\,\log 2c_k).\]
These choices ensure that
\[c_kq^{-\delta(2\sigma_0-1)}\le\frac{1}{2},\] 
and hence imply that
\[J(\sigma_0)\ll q^{4\delta(2\sigma_0-1)}\phi(q)(\log q)^{k^2-1}
\ll\phi(q)(\log q)^{k^2-1}.\]
Finally, we may apply Lemma \ref{Jconvlem2} to deduce the following
\begin{lemma}\label{Jest}
With $\sigma_0$ as above we have
\[J(\gamma)\ll\phi(q)(\log q)^{k^2-1}\]
uniformly for $1-\sigma_0\le \gamma\le \sigma_0$.
\end{lemma}

All that remains is to bound $M_k(q)$ from above, using averages of $J(\gamma)$.
Since $|L(s,\chi)|^{2k}$ is subharmonic we have
\[|L(\tfrac12,\chi)|^{2k}\le\frac{1}{2\pi}\int_0^{2\pi}
|L(\tfrac12+re^{i\theta},\chi)|^{2k}d\theta.\]
We now multiply by $r$ and integrate for $0\le r\le R$ to show that
\[|L(\tfrac12,\chi)|^{2k}\le\frac{1}{{\rm Meas}(D)}\int_D
|L(\tfrac12+z,\chi)|^{2k}dA,\]
where $D=D(0,R)$ is the disc of radius $R$ about the origin, and $dA$
is the measure of area.  We take 
\[R=\frac{\min\left(\kappa\,,\,\delta^{-1}\right)}{\log q},\]
so that if $z\in D$ then $1-\sigma_0\le\Re(\tfrac12+z)\le\sigma_0$ and
$|W(\tfrac12+z)|\gg 1$.  It follows that
\[\int_D |L(\tfrac12+z,\chi)|^{2k}dA\ll 
\int_{1-\sigma_0}^{\sigma_0}J(\gamma,\chi)d\gamma\]
whence
\[M_k(q)\ll
\frac{1}{{\rm Meas}(D)}\int_{1-\sigma_0}^{\sigma_0}J(\gamma)d\gamma.\]
Since ${\rm Meas}(D)\gg (\log q)^{-2}$ we now deduce from Lemma
\ref{Jest} that 
\[M_k(q)\ll\phi(q)(\log q)^{k^2},\]
as required.

\bigskip
\bigskip

Mathematical Institute,

24--29, St. Giles',

Oxford

OX1 3LB

UK
\bigskip

{\tt rhb@maths.ox.ac.uk}

\end{document}